\documentclass{article}

\newtheorem{fed}{\textbf{Definition}}[section]
\newtheorem{thm}[fed]{\textbf{Theorem}}
\newtheorem{lemma}[fed]{\textbf{Lemma}}

\usepackage{amssymb,bbm,graphicx,epsfig,psfrag,epic,eepic,latexsym}
\usepackage{amsmath}
\usepackage{mathrsfs}
\usepackage{comment}
\usepackage{hyperref}
\usepackage{color}

\newcommand{\Z}{\mathbb{Z}}
\newcommand{\R}{\mathbb{R}}

\newcommand{\Address}{{
  \bigskip
  \footnotesize 
  \textsc{Lehrstuhl für Analysis und Geometrie, Universität Augsburg, D-86135 Augsburg, Germany}\\
  \textit{E-mail addresses:} \href{mailto:filipvbrocic@gmail.com}{\tt filipvbrocic@gmail.com} \hspace{15mm} \href{mailto:urs.frauenfelder@uni-a.de}{\tt urs.frauenfelder@uni-a.de}
}}

\begin{document}
	\title{Wrapped Floer homology and the circular restricted three-body problem}
	\author{Filip Bro\'ci\'c and Urs Frauenfelder}
	\maketitle
	
	\begin{abstract}
		Using the wrapped Floer homology, we prove the existence of consecutive collisions at the primaries in the circular restricted three-body problem. We also prove the existence of a symmetric periodic orbit. These existence results are obtained for energy hypersurfaces slightly above the first critical value.
	\end{abstract}
	
	\section{Introduction}
	
	The purpose of this article is to provide a different viewpoint on a recent theorem by Kang and Ruck \cite{KR25} on the existence of consecutive collisions
	in the restricted three-body problem slightly above the first critical value. While the proof of Kang and Ruck uses explicit computations in equivariant Rabinowitz Floer homology, our proof is based on surgery techniques in wrapped Floer homology due to Irie \cite{Ir13} and the first author of this paper \cite{Br25}. In \cite{AFKP12}, it was shown that the planar restricted three-body problem below and slightly above the first critical value admits a contact structure. 
	After regularizing collisions for each energy value below the first critical value, the energy hypersurface can be interpreted as a fiberwise starshaped hypersurface in the cotangent bundle of $S^2$. In the regularization, the role of position and momentum was interchanged, and therefore the sphere corresponds to the momenta, and the fiber $T^*_\infty S^2$ over the north pole corresponds to collisions. In particular, by a result of Abbondandolo and Schwarz \cite{ASch06}, the wrapped Floer homology of this fiber is isomorphic to the homology of the based loop space of the sphere
	$$HW_*(T_{\infty} S^2) \cong H_*(\Omega(S^2)).$$
	In \cite{FZ19}, this fact was used to prove the existence of a periodic collisional orbit or infinitely many consecutive collisions. As discovered by Ruck \cite{Ru24}, finer information on consecutive collisions can be obtained by taking into account that the restricted three-body problem is invariant under an antisymplectic involution $\rho$.  
	In the original coordinates, the existence of this antisymplectic involution is based on the fact that the restricted three-body problem is invariant under reflection at the axis on which the primaries lie. For the definition of $\rho$ see \S\ref{sec:hamiltonian}. 
	\\ \\
	The fixed point set $\mathrm{Fix}(\rho)$ of an antisymplectic involution is a Lagrangian submanifold. In the regularized coordinates, the fixed point set corresponds to the conormal bundle of a meridian in $T^*S^2$. By a result of Abbondandolo, Portaluri, and Schwarz \cite{APS08}, it follows that the wrapped Floer homologies are given by
	\begin{equation}\label{colfix}
	HW_*(T_{\infty} S^2,\mathrm{Fix}(\rho)) \cong H_*(\mathcal{P}_{\infty, S^1})
	\end{equation}
	respectively
	\begin{equation}\label{fix}
	HW_*(\mathrm{Fix}(\rho)) \cong H_*(\mathcal{P}_{S^1, S^1}).
	\end{equation}
	Here $\mathcal{P}_{N,Q}$ denotes the space of paths in the ambient manifold, starting on $N$ and ending on $Q$. Since the meridian $S^1 \subset S^2$ is contractible in $S^2$ we have $\mathcal{P}_{\infty, S^1} \simeq \Omega(S^2) \times S^1$ and $\mathcal{P}_{S^1,S^1} \simeq \Omega(S^2) \times S^1 \times S^1$.
	\\ \\
Generators of $HW_*(T_{\infty} S^2,\mathrm{Fix}(\rho))$ are chords starting in a collision and ending in a fixed point of $\rho$. Since the collision Lagrangian $T_\infty^*S^2$ is invariant under the antisymplectic involution $\rho$, it follows that with the help of the antisymplectic involution, each such chord can be doubled to get a consecutive collision which has the additional property that it is invariant under simultaneously reversing time and applying the antisymplectic involution. These special types of consecutive collision orbits are referred by Ruck \cite{Ru24} as \emph{symmetric consecutive collisions}. In particular, it follows from \cite{Ru24}
	that there are infinitely many symmetric consecutive collisions or a periodic symmetric consecutive collisions for energies below the first critical value. By using the fact that the Hamiltonian of the restricted three-body problem is analytic in the mass parameter, an argument by Kang and Ruck \cite{KR25} even shows that for generic mass parameters, there have to exist at least two consecutive collisions. 
	\\ \\
	From (\ref{fix}) existence of chords from the fixed point set of the Hamiltonian involution to itself is obtained. Such chords can be doubled as well and give rise to symmetric periodic orbits. 
	\\ \\
	Slighly above the first critical value, energy hypersurfaces of the restricted three-body problem are still of contact type, and passing the critical value corresponds to contact connected sum \cite{AFKP12}. Above the first critical value, collisions with both primaries, the earth and the moon, are possible. We refer to $L_e$ as the Lagrangian corresponding to collisions with the earth and to $L_m$ as the Lagrangian corresponding to collisions with the moon. We still have the antisymplectic involution $\rho$ and its fixed point set $\mathrm{Fix}(\rho)$. Using surgery in wrapped Floer homology as in \cite{Br25, Ir13} we compute the corresponding wrapped Floer homologies. Namely
	
	\begin{thm}
	Slightly above the first critical value, the wrapped Floer homologies of the three Lagrangians $L_e$, $L_m$, and $\mathrm{Fix}(\rho)$ are given by
	\begin{equation}\label{colfix+}
	HW_*(L_e,\mathrm{Fix}(\rho)) \cong HW_*(L_m,\mathrm{Fix}(\rho)) \cong H_*(\Omega(S^2)\times S^1 )
	\end{equation}	
	\begin{equation}\label{em}
	HW_*(L_e,L_m)\cong HW_*(L_m,L_e) \cong \{0\}	
	\end{equation}
	\begin{equation}\label{col}
	HW_*(L_e) \cong HW_*(L_m) \cong HW_*(\Omega(S^2))	
	\end{equation}
	\begin{equation}\label{fix+}
	HW_*(\mathrm{Fix}(\rho)) \cong H_*(\Omega(S^2) \times S^1 \times S^1) \oplus H_*(\Omega(S^2) \times S^1 \times S^1).
	\end{equation}	
	\end{thm}
	From (\ref{col}) we obtain an alternative proof of the theorem of Kang and Ruck of the existence of infinitely many consecutive collisions or a periodic consecutive collision where both collisions are at the moon, respectively, the earth.  From (\ref{fix+}) it follows that symmetric periodic orbits also have to exist slightly above the first critical value. From (\ref{colfix+}) we get that there are either infinitely many symmetric consecutive collisions or a symmetric periodic consecutive collision. Generators of $HW_*(L_e,L_m)$ are consecutive collisions where the first collision is with the earth and the other with the moon, and similarly for $HW_*(L_m,L_e)$ with the roles of earth and moon interchanged. Since this homology vanishes, one cannot directly deduce the existence of consecutive collisions starting at one of the primaries and ending at the other one. However, the computation of this homology gives rise to a forcing result - if there is one nondegenerate consecutive collision from the earth to the moon, there has to be a second one, and similarly for consecutive collisions starting at the moon and ending at the earth. It is an intriguing question if it might be possible to show that the chain complexes are nevertheless nontrivial, although the homology vanishes.  This would require a deeper understanding of the Fukaya category of the restricted three-body problem slightly above the first critical value. We expect that the careful study of Lagrangian subcritical surgery of the first author of this paper \cite{Br25}, which led to the current note, can as well give insight into the underlying Fukaya categories.
	\subsection{On the proof}
	Even though the strategy is clear, here we include how to collect all the pieces. The Lagrangians $L_e$ and $L_m$ are not affected by the surgery. However, the Lagrangian $\mathrm{Fix}(\rho)$ before the surgery has two connected components $\nu^* S^1_e \sqcup \nu^*S^1_m$, each one being inside the corresponding component of the filling $T^*S^2_e \sqcup T^*S^2_m$ of the regularized energy hypersurfaces $\overline{\Sigma}_c^e \sqcup \overline{\Sigma}_c^m$. After the surgery, $\mathrm{Fix}(\rho)$ is diffeomorphic to the pair of pants surface. This is the case since both components of $\mathrm{Fix}(\rho)$ before the surgery contain the points in  $\overline{\Sigma}_c^e \sqcup \overline{\Sigma}_c^m$ where the surgery is taking place.
\\ \\	
Now, denote by $W$ the filling of $\overline{\Sigma}^{e,m}_c$. It follows from Theorem \ref{thm:invariance} that we have
	$$
	HW_*(L_e, \mathrm{Fix}(\rho); W) \cong HW_*(L_e, \nu^*S^1_e \sqcup \nu^* S^1_m ; T^*S^2_e \sqcup T^*S^2_m).
	$$
	Since before the surgery $L_e  = T^*_\infty S^2_e$ and $ \nu^* S^1_m \subset T^*S^2_m$ do not interact we get
	$$
	HW_*(L_e, \mathrm{Fix}(\rho); W) \cong HW_*(T^*_{\infty} S^2, \nu^* S^1) \cong H_* (\Omega(S^2) \times S^1).
	$$ 
	In the same fashion, one gets the result for $HW_*(L_m, \mathrm{Fix}(\rho))$, which proves part (\ref{colfix+}). For part (\ref{em}), we get the result since $L_e$ and $L_m$ do not interact before the surgery, hence by the invariance theorem we have $HW_*(L_e, L_m;W) \cong HW_*(L_m, L_e;W) \cong 0$. Since $L_e$, and $L_m$ are not affected by the surgery we get from Theorem \ref{thm:invariance}:
	$$
	HW_*(L_e; W) \cong HW_*(L_e; T^*S^2_e \sqcup T^* S^2_m) = HW_*(T^*_{\infty} S^2_e; T^*S^2_e) \cong H_*(\Omega(S^2)),
	$$
	and by the same argument we get the result for $HW_*(L_m)$, concluding part (\ref{col}). For $HW_*(\mathrm{Fix}(\rho); W)$ we get
	$$
	\begin{aligned}
	HW_*(\mathrm{Fix}(\rho); W) &\cong HW_*(\nu^*S^1_e \sqcup \nu^*S^1_m; T^*S^2_e \sqcup T^*S^2_m)\\ &= HW_*(\nu^*S^1_e; T^*S^2_e) \oplus HW_*(\nu^*S^1_m; T^*S^2_m),
	\end{aligned}
$$
combining this with (\ref{fix}) leads to part (\ref{fix+}). 

\subsection{Organization of the paper}
In \S\ref{sec:wrapped_floer}, we recall the definition of the wrapped Floer homology. In \S\ref{sec:handle_attach}, we cover the handle attachment and the invariance theorem for wrapped Floer homology under subcritical handle attachment. In \S\ref{sec:R3BP} we give a background on the restricted three-body problem, and in \S\ref{sec:contact_str_ham}, we describe what happens with the bounded component of the energy hypersurface above the first critical value.

\section*{Acknowledgmenets}
This project was supported by the Deutsche Forschungsgemeinschaft (DFG, German Research Foundation) – 517480394.
	
	\section{Wrapped Floer homology}\label{sec:wrapped_floer}
	The wrapped Floer homology is a Lagrangian counterpart of the symplectic homology. It was introduced in \cite{ASch06} for the case of a fiber $T_q^* N$ in the cotangent bundle $T^*N$, extended to the case of conormal bundle $\nu^*Q$ of a submanifold $Q \subset N$ in \cite{APS08} and generalized in \cite{AS10} to the class of exact \textit{cylindrical} Lagrangians in Liouville manifolds. Here we give a brief overview of the definition of the wrapped Floer homology; for more details, see \cite{Ab12, Ri13, Ir13, BCS24, Br25}.
\\ \\ 	
A triple $(M, \lambda, X)$ is a \textit{Liouville domain} with a contact boundary $\partial M$ if $(M, \omega:=d \lambda)$ is a compact symplectic manifold with boundary $\partial M$, and if the vector field $X$ given by $i_X \omega = \lambda$ is positively transverse to $\partial M$. Transversality of $X$ to $\partial M$ insures that the form $\alpha:= \lambda \vert_{\partial M}$ is contact. The Liouville vector field $X$ can be used to complete $(M, \lambda, X)$ to a Liouville domain by attaching the positive symplectization $(\partial M \times [1,\infty), d(r \alpha))$ of $(\partial M, \alpha)$ to the boundary of $(M, d \lambda)$. We denote the completion of a Liouville domain by $(\widehat{M}, d\widehat{\lambda}, \widehat{X})$. Note that $\widehat{X} = \partial_r$ on $\partial M \times [1, \infty)$. We will consider time-dependent, $\widehat{\omega}$-compatible almost complex structures $J$ which are time independent in the symplectization and satisfy 
	\begin{equation}\label{eq:sft_acs}
		\widehat{\lambda} \circ J = d r \text{ on } \{r\geq1\}.
	\end{equation}
A Lagrangian $\widehat{L} \subset \widehat{M}$ is called \textit{cylindrical} if for all $p = (x,r) \in \partial M \times [1, \infty)$ we have $\widehat{X}(p) \in T_p \widehat{L}$. Equivalently, there is a Legendrian $\Lambda \subset \partial M$ such that $\widehat{L} \cap \partial M \times [1, \infty) = \Lambda \times [1, \infty)$. For our purposes, it is enough to assume that the Liouville form $\widehat{\lambda}$ vanishes along Lagrangians, in general, it is assumed that $\widehat{\lambda}$ is exact when restricted on $\widehat{L}$ meaning that there is a function $f_{\widehat{L}}: \widehat{L} \to \R$ such that $\widehat{\lambda} \vert_{\widehat{L}} = d f_{\widehat{L}}$. Since $\widehat{L}$ is cylindrical $f_{\widehat{L}}$ is locally constant. Additionally, we assume that the relative Chern class 
	\begin{equation}\label{eq:rel_chern}
		2c_1(M,L) \in H^2(M, L)  \text{ vanishes}.
	\end{equation} 
	This assumption ensures that we can have a globally defined $\Z$ grading on the Floer groups. We could also assume that the Maslov class $\mu_L \in H^1(L)$ vanishes, which implies $2c_1(M,L)=0$, and is sufficient for our applications.
\\ \\	
Let $R_{\alpha}$ be the Reeb vector field for $(\partial M, \alpha)$. Given a Legendrian $\Lambda$  define the \textit{action spectrum} of $(\Lambda, \alpha)$ :
	$$
	\mathcal{A}(\Lambda, \alpha) := \left\{ \int \gamma^* \alpha \mid \gamma \text{ is a Reeb chord with boundary on } \Lambda\right\}.
	$$
	Analogously, define the action spectrum for two Legendrians  $\Lambda_0, \Lambda_1$ to be
	$$
	\mathcal{A}(\Lambda_0, \Lambda_1, \alpha) := \left\{ \int \gamma^* \alpha \mid \gamma \text{ is a Reeb chord from } \Lambda_0 \text{ to } \Lambda_1\right\}.
	$$
	Given a Hamiltonian function $H:[0,1] \times \widehat{M} \to \R$, the Hamiltonian vector field $X_{H_t}$ is defined by $i_{X_{H_t}} \widehat{\omega} = -dH_t$. A Hamiltonian is called \textit{admissible} for $\widehat{L}$ (resp. $(\widehat{L}_0, \widehat{L}_1)$) if 
	$$
	H(x,r) = a r + b, \text{ on } \partial M \times [1, \infty),
	$$
	for constants $a>0$,b and $a \notin \mathcal{A}(\Lambda, \alpha)$ (resp. $a \notin \mathcal{A}(\Lambda_0, \Lambda_1, \alpha)$). We set $\mathcal{H}$ to be the set of all admissible Hamiltonians. We call the number $a>0$ the slope of $H_t$. Hamiltonian $H_t$ is \textit{non-degenerate} if $\varphi^1_{H_t}(L) \pitchfork L$ (resp. $\varphi^1_{H_t}(L_0) \pitchfork L_1$), here $\varphi^t_{H_t}$ is the flow of $X_{H_t}$. The assumption that the slope $a$ is not in the action spectrum follows from the non-degeneracy of a Hamiltonian which has slope $a$ at infinity. Let $\mathcal{H}^{\mathrm{x}} \subset \mathcal{H}$ be the set of all non-degenerate admissible Hamiltonians.
\\ \\
The action functional on the space of paths $$\mathcal{P}_L:=\{x:[0,1] \to \widehat{M} \mid x(0), x(1) \in \widehat{L}\}$$ is defined by
	\begin{equation}\label{eq:act_fun}
		\mathcal{A}_{H_t}(x) = \int x^* \widehat{\lambda} - \int_0^1 H_t(x(t)) dt.
	\end{equation}
	The critical points of $\mathcal{A}_{H_t}$ are Hamiltonian paths $x:[0,1] \to \widehat{M}$ with $x(0), x(1) \in L$. Here, we have used that $\widehat{\lambda} \vert_L =0$, for a general exact Lagrangian, the action functional should involve the term $f_L(x(1)) - f_L(x(0))$. We denote the set of critical points $\mathrm{Crit}(\mathcal{A}_{H_t})$. Analogously one defines the action functional on the space of paths $\mathcal{P}_{L_0, L_1}$ with endpoints on $L_0$ and $L_1$ whose critical points are Hamiltonian paths $x:[0,1] \to \widehat{M}$ with $x(i) \in L_i$ for $i \in \{0,1\}$. Since the slope is not in the action spectrum, the images of all elements in $\mathrm{Crit}A_{H_t}$ are contained in the compact part $M \subset \widehat{M}$.
\\ \\	
The Floer chain group for an admissible Hamiltonian $H$ is defined by
	\begin{equation}\label{eq:chain_complex}
		CF_k(L,H,J) = \bigoplus_{\substack{x \in \mathrm{Crit}A_{H_t}, \\ \mu(x)=k}} \Z_2 \langle x \rangle,
	\end{equation}
	here $\mu(x)$ is the Maslov index. We normalize $\mu$ so that if $x$ is a critical point of $H$ which is a $C^2$ small extension of a $C^2$ small Morse function on $L$ then $\mu(x) =n-m_f(x)$ where $m_f(x)$ is the Morse index of $x$ as a critical point of $f$. 
\\ \\	
Given $x_-, x_+ \in \mathrm{Crit}\mathcal{A}_{H_t}$ the moduli space $\mathcal{M}(x_-,x_+,H,J)$ is the set of $u:\R \times [0,1] \to \widehat{M}$ which satisfy the \textit{Floer equation}:
	
	\begin{equation}\label{eq:floer_eq}
		\partial_s u + J_t (\partial_t u - X_{H_t}) = 0,
	\end{equation}
	with the boundary conditions $u(s,0), u(s,1) \in \widehat{L}$, and with asymptotic conditions $\lim\limits_{s \to \pm \infty} u(s,t) = x_{\pm}(t)$. For a generic $H_t$ the moduli space $\mathcal{M}(x_-,x_+,H,J)$ is a smooth manifold of dimension $\mu(x_+) - \mu(x_-)$. There is an $\R$ action on $\mathcal{M}(x_-,x_+,H,J)$ by translation in $s$ direction. Denote $$\Bar{\mathcal{M}}(x_-, x_+) = \mathcal{M}(x_-,x_+,H,J) / \R.$$ When the index difference $\mu(x_+) - \mu(x_-)$ is $1$, standard compactness results  imply that $\Bar{\mathcal{M}}(x_-, x_+)$ is a finite number of points.
\\ \\	
The differential $d:CF_k(L,H,J) \to CF_{k-1}(L,H,J)$ is defined on generators $y \in \mathrm{Crit}\mathcal{A}_{H_t}$  by
	\begin{equation}\label{eq:dif}
		d y = \sum_{\substack{x,\\ \mu(x) = k-1}} \#_2 \Bar{\mathcal{M}}(x, y) x,
	\end{equation}
	and extended to $CF_k(L,H,J)$ by linearity. By the standard gluing and compactness arguments, we have that $d^2=0$. We denote $HF_k(L,H,J)$ the homology of the complex $(CF_k(L,H,J), d)$
\\ \\	
Fix two admissible Hamiltonians $H_-$ and $H_+$ with the slopes $a_-$ and $a_+$. If $a_+ \leq a_-$, one can define a continuation map 
	\begin{equation}\label{eq:cont_maps}
		\Phi_{H_s}:CF_k(L, H_+,J_+) \to CF_k(L,H_-,J_-)
	\end{equation}	
	Let $\mathcal{M}(x_-,x_+,H_s,J_s)$ be the moduli space  of solutions
	\begin{equation}\label{eq:cont_eq}
		\partial_s u + J_{s,t} (\partial_t u - X_{H_{s,t}}) = 0,
	\end{equation}
	with the same boundary conditions and analogous asymptotic conditions as in (\ref{eq:floer_eq}). Here $H_{s,t}$ is a homotopy between $H_-$ and $H_+$ satisfying $H_s = H_-$ for $s \leq -s_0$, $H_s = H_+$ for $s \geq s_0$ and $\partial_s H_s \leq 0$. The last condition ensures that we can apply the maximum principle, which is needed for the compactness of the moduli space $\mathcal{M}(x_-,x_+,H_s,J_s)$. Domain dependent almost complex structure $J_{s,t}$ is a homotopy between $J_-$ and $J_+$. For a generic choice homotopy $H_s$, the moduli space is a smooth manifold of dimension $\mu(x_+) - \mu(x_-)$. The map $\Phi_{H_s}:CF_k(L, H_+,J_+) \to CF_k(L,H_-,J_-)$ is defined on generators by
	\begin{equation}\label{eq:cont_def}
		\Phi_{H_s}(y) = \sum_{\substack{x,\\ \mu(x) = \mu(y)}} \#_2 \mathcal{M}(x, y, H_s, J_s) x,
	\end{equation}
	Note that since $H_s$ depends on $s$, solutions are not translation invariant, so we do not have an $\R$ action. By the standard gluing and compactness arguments, one shows that $\Phi_{H_s}$ is a chain map, and we denote again with $\Phi_{H_s}$ the map induced on the homology. By considering homotopies of homotopies, one can show that the map $\Phi_{H_s}$ on homology does not depend on the choice of $H_s$.
\\ \\	
Analogously one defines the Floer groups $CF_k(L_0, L_1, H, J)$ and continuation maps $\Phi_{H_s}:CF_k(L_0, L_1, H_+,J_+) \to CF_k(L_0, L_1,H_-,J_-)$, with appropriate change of the boundary conditions in (\ref{eq:floer_eq}) and in (\ref{eq:cont_eq}). 
\\ \\ 	
On the space of admissible Hamiltonians $\mathcal{H}$, there is a partial order 
	$$
	H_1 \preccurlyeq H_2 \text{ if and only if } a_1 \leq a_2,
	$$
	where $a_i$ is the slope of $H_i$. The continuation map $\Phi_s: HF(L,H_1) \to HF(L,H_2)$ is well defined if $H_1 \preccurlyeq H_2$. The wrapped Floer homology of $L$ is defined by
	\begin{equation}\label{eq:wrapped_def}
		HW_*(L;M) := {\lim_{\substack{\longrightarrow \\ H \in \mathcal{H}^{\mathrm{x}}}}} HF_*(L,H).
	\end{equation}
	Similarly, one defines the wrapped Floer homology of the pair $(L_0, L_1)$:
	\begin{equation}\label{eq:wrapped_def}
		HW(L_0, L_1;M)_* := {\lim_{\substack{\longrightarrow \\ H \in \mathcal{H}^{\mathrm{x}}}}} HF_*(L_0, L_1,H).
	\end{equation}
	If we choose a sequence of Hamiltonians $H_i \in \mathcal{H}^{\mathrm{x}}$ such that the slopes $a_i$ form an increasing unbounded sequence, then $H_i$ is a cofinal sequence in $\mathcal{H}^{\mathrm{x}}$ and we have
	$$
	HW(L) \cong \lim_{\substack{\longrightarrow\\i}} HF(L, H_i).
	$$
	The same argument holds for $HW(L_0, L_1)$.
\\ \\
An alternative approach to define wrapped Floer homology groups is to take a single Hamiltonian $H$ which is of the form $\frac{1}{2}r^2$ on $\partial M \times [1,\infty)$. One can show that $HF(L, H) \cong HW(L)$, for this comparison see \cite{Ri13, BCS24}. The approach with the slopes is more convenient for us for the invariance of the wrapped Floer homology group under the subcritical surgeries, as we will see in section \ref{sec:handle_attach}. The approach with the quadratic Hamiltonian was used in \cite{ASch06, APS08}.
\\ \\	
For us, the important case is when $\widehat{M}$ is the cotangent bundle $T^* N$ of a closed smooth manifold $N$ with the standard symplectic structure, and $\widehat{L_i}$ are conormal bundles $\nu^* Q_i$ of closed submanifolds $Q_i \subset N$. In particular, we will consider the case $N = S^2$ and $Q_i$ are either the north pole $\{N\} \in S^2$ or a great circle $S^1 \subset S^2$ passing through the north and south poles.
	Let $\mathcal{P}_{Q_0, Q_1}$ be the space of $W^{1,2}$ paths $x:[0,1] \to N$ with $x(i) \in Q_i$ for $i \in \{0,1\}$. In \cite{APS08} they showed the following
	\begin{thm}\label{thm:conormal_iso}
		$$
		HW_*(\nu^* Q_0, \nu^* Q_1) \cong H_*(\mathcal{P}_{Q_0, Q_1}).
		$$
	\end{thm}
	The special case when both $Q_i$ are points $q_i$ in $N$ was established in \cite{ASch06}, then $\mathcal{P}_{q_0, q_1}$ is homotopy equivalent to the based loop space $\Omega_q(N)$. Theorem \ref{thm:conormal_iso} holds with arbitrary coefficients; for our purpose, it is enough to consider the $\Z_2$ coefficients. In this case, we can avoid additional complications with the coherent orientations of the moduli spaces.
	
	\subsection{Handle attachment}\label{sec:handle_attach}
	
	In \cite{Ir13}, it was shown that the wrapped Floer homology is invariant under the subcritical handle attachment. This is motivated by \cite{Ci02}, where it was shown that the same holds for symplectic homology. In \cite{Fa16-phd, Fa20} Fauck discovered and corrected a small gap in the choice of the cofinal family of Hamiltonians in \cite{Ci02}. This gap was carried over in \cite{Ir13}. Following ideas from \cite{Fa16-phd, Fa20}, this was corrected in an expository article \cite{Br25}.
\\ \\	
Recall the construction of the contact surgery from \cite{We91}. Let $S \cong S^{k-1}$ be an isotropic sphere of a contact manifold $(Y, \xi )$ with a trivial conformal symplectic normal bundle $CSN(S):= (TS)^{\perp_{\omega}} / TS$. Here $\omega$ is the symplectic structure on $\xi$ uniquely determined up to a conformal factor. One obtains a new manifold $Y'$ from $Y$ by the surgery along $S$. This new manifold can be equipped with a contact structure, and the trace of the surgery $W$ is a Liouville cobordism between $Y$ and $Y'$.
\\ \\	
To describe this more precisely, we need to introduce a standard handle. On $(\R^{2n}, \omega_{st})$ consider a Liouville form given by
	
	\begin{equation*}
		\phi = \sum_{i=1}^{k} (2x_i dy_i + y_i dx_i)+ \frac{1}{2}\sum_{i=k+1}^n (x_i dy_i - y_i dx_i),
	\end{equation*}
	the corresponding Liouville vector field is:
	$$
	X= \sum_{i=1}^{k} (2x_i \partial_{x_i} - y_i \partial_{y_i})+ \frac{1}{2}\sum_{i=k+1}^n (x_i \partial_{x_i} + y_i \partial_{y_i}).
	$$
	The handle will be a region between two contact type hypersurfaces, that are transverse to $X$. These hypersurfaces are determined by functions $\phi$ and $\psi$. First, let us define: 
	$$
	\phi = \frac{1}{2}\sum_{i=1}^{k} (2x_i^2 + y_i^2)+ \frac{1}{4}\sum_{i=k+1}^n (x_i^2 + y_i^2),
	$$
	and $\Sigma_- = \{ \phi = -1\}$. A neighborhood $U$ of the isotropic sphere $S = \{x=z=0, y=1\}$ is identified with a neighborhood of $S \subset Y$. Hence, we will choose so that $\Sigma_-$ coincides with $\Sigma_+ = \{\psi = -1\}$ outside of $U$. Following \cite{Fa20, Br25}, for $\epsilon$ we set
	$$
	g(t)= \begin{cases}
		\frac{1}{1+2\epsilon} t, &t \leq 1,\\
		1, &t\geq 1+ 3\epsilon,
	\end{cases}
	$$ 
	so that $0 \geq g' \geq 1/(1+2 \epsilon)$. Define functions $x,y,z : \R^{2n} \to \R$ by
	$$
	x= \sum_{i=1}^k x_i^2, \hspace{5mm} y = \frac{1}{2} \sum_{i=1}^k  y_i^2, \hspace{5mm} z = \frac{1}{4} \sum_{i=k+1}^{n} \left(x_i^2 + y_i^2 \right).
	$$
	The function $\phi$ takes the form $\phi = x-y+z$. Given $\delta >0$ set
	$$
	\psi_{\delta} = x-y+z - (1+\epsilon) + (1+\epsilon) g(y + (x+z)/\delta).
	$$
	The parameter $\delta$ is used to ensure that $\Sigma_+$ defined by $\psi:=\psi_{\delta}$ and $\Sigma_-$ coincide outside of a small neighborhood of $S$. The handle is defined by 
	$$
	H_k^{2n}:= \{\phi \geq -1\} \cap \{\psi_{\delta} \leq -1\}.
	$$
	Using the Liouville vector field $X$ on the handle, we see that the outcome of the surgery carries a contact structure, i.e., we have the following
	\begin{thm}[\cite{We91}]\label{thm:handle_attach}
		Let $(Y,\xi)$ be a closed manifold with co-oriented contact structure, and let $S \subset Y$ be an isotropic sphere with trivialized conformal symplectic normal bundle $CSN(S)$. The outcome $Y'$ of the surgery of $Y$ along $S$ carries a contact structure. Furthermore, there exists a Liouville cobordism from $Y$ to $Y'$.
	\end{thm}
	It follows from \cite[Lemma 2.6]{Ci02} that if an isotropic sphere $S$ is contained in the Legendrian $\Lambda \subset Y$, one gets an exact Lagrangian cobordism from $\Lambda$ to $\Lambda'$, where $\Lambda'$ is obtained by surgery on $\Lambda$, along $S$, by attaching $D^{k}\times S^{n-k-1} \times \{0\} \subset $ and the Lagrangian cobordism is obtained by attaching the imaginary part of the handle $H_k^{2n}$ given by $H_k^{n}:= \{(0,...,0, y_1,...,y_n) \in H_k^{2n}\}$. Now we are in a position to state the invariance theorem. 
	
	\begin{thm}[\cite{Ir13, Br25}]\label{thm:invariance}
		Given two exact, graded, cylindrical Lagrangians $L_i \subset M$, we have 
		$$
		HW_*(L_0, L_1; M) \cong HW_*(L_0, L_1 \cup_{S} H_k^n; M \cup_S H_k^{2n}).
		$$ 
		In the case when $L=L_0 = L_1$ we have $HW_*(L;M) \cong HW_*(L \cup H_k^n; M \cup H_k^{2n})$.
	\end{thm}
	One can see from the proof that it is also possible to allow any Lagrangian plane in the handle that passes through origin, so, it is possible to take $L_1 = \emptyset$, or, $L_1$ can intersect $S$ along some sphere of lower dimension in such a way that the glued Lagrangian coincides with the Lagrangian plane in the handle, see \cite[Remark 1.1]{Br25}.

	\section{The restricted three-body problem}\label{sec:R3BP}
	
	\subsection{The Hamiltonian}\label{sec:hamiltonian}
	
	The restricted three-body problem describes the motion of a massless particle attracted according to Newton's law of gravitation by two massive bodies, the primaries. We refer to the two primaries as the earth and the moon and to the massless body as the satellite. In the circular case, one assumes that in the inertial system,
	the earth and the moon move on a circle around their common center of mass. We further consider the planar case in which the satellite moves in the eccliptic, i.e., the same plane as the earth and the moon. Since the earth and the moon are moving, the Hamiltonian for the satellite depends periodically on time in the inertial system. For that reason, one considers the system in rotating coordinates in which both the earth and the moon are at rest. In such a rotating system, the Hamiltonian does not depend on time anymore, i.e., it is autonomous and therefore preserved along the flow of its Hamiltonian vector field. We scale the total mass of the earth and the moon to one, and denote by  $\mu \in (0,1)$ the mass of the moon. The mass of the earth is then given by $1-\mu$. We assume without loss of generality in the following that $\mu \leq \tfrac{1}{2}$, otherwise we just interchange the roles of the moon and the earth. As our length unit, we choose the distance between the earth and the moon, and choose our coordinates such that the earth lies at $e=(-\mu,0) \in \mathbb{R}^2$
	and the moon at $m=(1-\mu,0) \in \mathbb{R}^2$. The center of mass lies at the origin. The phase space for the satellite is then the cotangent bundle
	$$T^* \big(\mathbb{R}^2 \setminus \{e,m\}\big)=\big(\mathbb{R}^2 \setminus \{e,m\}\}\big) \times \mathbb{R}^2$$
	and the Hamiltonian for the satellite becomes
	$$H \colon T^* (\mathbb{R}^2 \setminus \{e,m\}) \to \mathbb{R}, \quad (q,p) \mapsto \frac{1}{2}||p||^2-\frac{\mu}{||q-m||}-\frac{1-\mu}{||q-e||}+p_1q_2-p_2q_1.$$
	The first term is the kinetic energy of the satellite, the second term is the gravitational potential of the moon, and the third term is the gravitation potential of the earth. A bit more surprising is the last term $p_1q_2-p_2q_1$. This last term is the angular momentum whose Hamiltonian vector field generates the rotation. Since we are considering the system in rotating coordinates and not in the inertial system, we have to add this term to the Hamiltonian. In particular, due to the addition of angular momentum, the Hamiltonian of the restricted three-body problem is not a mechanical Hamiltonian consisting just of kinetic and potential energy. However, we can complete the squares and rewrite the Hamiltonian as follows 
	$$H(q,p)=\frac{1}{2}\Big((p_1+q_2)^2+(p_2-q_1)^2\Big)-\frac{\mu}{||q-m||}-\frac{1-\mu}{||q-e||}-\frac{1}{2}||q||^2.$$
	The last three terms only depend on position $q$ and not on momentum $p$. We summarize them in the function
	$$U \colon \mathbb{R}^2 \setminus \{e,m\} \to \mathbb{R}, \quad q \mapsto -\frac{\mu}{||q-m||}-\frac{1-\mu}{||q-e||}-\frac{1}{2}||q||^2$$
	referred to as the \emph{effective potential}. The effective potential not just consists of the gravitational potentials of the earth and moon, but of an additional third term $-\tfrac{1}{2}||q||^2$ which gives rise to the centrifugal force experienced in a rotating system. With the help of the effective potential, we can rewrite the Hamiltonian as
	\begin{equation}\label{ham}
		H(q,p)=\frac{1}{2}\Big((p_1+q_2)^2+(p_2-q_1)^2\Big)+U(q).
	\end{equation}
	We see that the Hamiltonian of the restricted three-body problem is an example of a magnetic Hamiltonian, consisting of a potential and a twisted kinetic energy. The physical interpretation of the twist in the kinetic energy is the Coriolis force. In contrast to centrifugal force, the Coriolis force is velocity dependent, as the Lorentz force of a magnetic field. 
	\\ \\
	Another remarkable property of the Hamiltonian of the restricted three-body problem is that it is invariant under the antisymplectic involution 
	$$\rho \colon T^* \mathbb{R} \to T^* \mathbb{R} \quad (q_1,q_2,p_1,p_2) \to (q_1,-q_2,-p_1,p_2)$$
	which on configuration space $\mathbb{R}^2$ restricts to orthogonal reflection at the axis through earth and moon under which the effective potential is invariant, so that we get
	$$H \circ \rho=H.$$
	We refer to \cite{FVK18} for more details.
	
	\subsection{Lagrange points and Hill's regions}
	
	By (\ref{ham}) we see that the footpoint projection
	$$\pi \colon T^* \mathbb{R}^2 \to \mathbb{R}^2, \quad (q,p) \mapsto q$$
	from phase space to configuration space gives rise to a bijection between critical points
	$$\Pi:=\pi\big|_{\mathrm{crit}(H)} \colon \mathrm{crit}(H) \to \mathrm{crit}(U)$$
	with inverse 
	$$\Pi^{-1} \colon \mathrm{crit}(U) \to \mathrm{crit}(H), \quad (q_1,q_2) \mapsto (q_1,q_2,-q_2,q_1).$$
	There are five critical points of $U$ referred to as the five Lagrange points $\ell_1,\ldots,\ell_5$. The first three Lagrange points $\ell_1,\ell_2$, and
	$\ell_3$ lies on the x-axis, i.e., on the axis through the earth and the moon. The first Lagrange point $\ell_1$ lies between the earth and the moon, $\ell_2$ lies to the right of the moon, and $\ell_3$ lies to the left of the earth, so that they are arranged in the following
	order
	$$\ell_3<e<\ell_1<m<\ell_2.$$
	They are nondegenerate saddle points of the effective potential
	$U$ and therefore in view of (\ref{ham}) their preimages under
	$\Pi$ are critical points of Morse index one of the Hamiltonian 
	$H$. The Lagrange points $\ell_4$ and $\ell_5$ are nondegenerate maxima of the effective potential and therefore their preimages under $\Pi$ are critical points of Morse index two of $H$. Each of them forms together with the earth and the moon an equilateral triangle. In particular, reflection at the x-axis interchanges $\ell_4$ and $\ell_5$ in agreement with the fact that the effective potential is invariant under reflection at the x-axis. 
	\\ \\
	From (\ref{ham}) we infer that the critical values of the effective potential and the Hamiltonian coincide, i.e., we have
	\begin{equation}\label{action}
		U(\ell_i)=H(\Pi^{-1}(\ell_i)), \quad 1 \leq i \leq 5.
	\end{equation}
	In the case $\mu<\tfrac{1}{2}$, i.e., the moon is less heavy than the earth the critical values of $U$ and therefore in view of (\ref{action}) also the ones of $H$ satisfy
	$$U(\ell_1)<U(\ell_2)<U(\ell_3)<U(\ell_4)=U(\ell_5).$$
	In the case where $\mu=\tfrac{1}{2}$ the earth and the moon have the same weight. This leads to an additional symmetry since now the earth and the moon can be interchanged. In this case, the effective potential is not only invariant under reflection at the x-axis, but as well under reflection at the y-axis, and the Hamiltonian $H$ becomes invariant under the antisymplectic involution
	$$\sigma \colon T^* \mathbb{R}^2 \to T^* \mathbb{R}^2, \quad
	(q_1,q_2,p_1,p_2) \mapsto (-q_1,q_2,p_1,-p_2)$$
	which commutes with the antisymplectic involution $\rho$ and whose product $\rho \sigma=\sigma \rho$ coincides with the symplectic involution $(q,p) \mapsto (-q,-p)$ on $T^* \mathbb{R}^2$.
	In this case, $\ell_1=0$ coincides with the barycenter of the earth and the moon, and reflection at the y-axis interchanges the second and third Lagrange point, so that we have
	$$\ell_2=-\ell_3.$$
	In the case $\mu=\tfrac{1}{2}$, there are not four but only three critical values, and one has
	$$U(\ell_1)<U(\ell_2)=U(\ell_3)<U(\ell_4)=U(\ell_5).$$
	Since the Hamiltonian $H$ is autonomous, it is preserved under its Hamiltonian flow. In particular, for any $c \in \mathbb{R}$ the level set
	$$\Sigma_c=H^{-1}(c) \subset T^*(\mathbb{R}^2 \setminus \{e,m\})$$
	is invariant under the flow. If $c$ is a regular value of $H$, and therefore by (\ref{action}) of $U$, the level set $\Sigma_c$ is a three-dimensional hypersurface in the four-dimensional phase space
	referred to as an energy hypersurface. 
	To understand the topology of these energy hypersurfaces, we introduce the notion of \emph{Hill's region}. The Hill's region is just the shadow of the energy hypersurface on configuration space, namely
	$$\mathfrak{K}_c:=\pi(\Sigma_c) \subset \mathbb{R}^2 \setminus \{e,m\}.$$
	In view of (\ref{ham}) it coincides with a sublevel set of the effective potential
	$$\mathfrak{K}_c=\big\{q \in \mathbb{R}^2 \setminus \{e,m\}: U(q) \leq c\big\}.$$
	The boundary of the Hill's region $\partial \mathfrak{K}_c$ has a special physical significance. These are the \emph{zero-velocity curves} where the satellite stops for an instance but is immediately accelerated again. 
	\\ \\
	If the energy $c$ lies below the first critical value $U(\ell_1)$ the Hill's region $\mathfrak{K}_c$ has three connected components
	$$\mathfrak{K}_c=\mathfrak{K}_c^e \cup \mathfrak{K}_c^m \cup \mathfrak{K}_c^u.$$
	The first two of them are bounded, the earth $e$ lies in the closure of 
	$\mathfrak{K}_c^e$, and the moon $m$ lies in the closure of $\mathfrak{K}_c^m$, i.e., the satellite either is close to the earth, close to the moon or is a comet and lies in the unbounded component $\mathfrak{K}_c^u$. Below the first critical value the satellite cannot travel from the earth to the moon, since the earth and the moon lie in the closures of different connected components of the Hill's region. This changes if the energy lies above the first critical value. In this case, a neck opens at the first Lagrange point, the two regions $\mathfrak{K}_c^e$ and $\mathfrak{K}_c^m$ get connected and there is no obvious obstruction anymore for the satellite to travel from the earth to the moon.

	\subsection{Regularization and contact structure below the first critical value}\label{sec:below}
	
	For $c<U(\ell_1)$ we abbreviate
	$$\Sigma_c^e:=\{(q,p) \in \Sigma_c: q \in \mathfrak{K}_c^e\}, \quad 
	\Sigma_c^m:=\{(q,p) \in \Sigma_c: q \in \mathfrak{K}_c^m\}$$
	the components of the energy hypersurface around the earth and the moon. 
	Although the Hill's regions $\mathfrak{K}_c^e$ and $\mathfrak{K}_c^m$ are bounded the components $\Sigma_c^e$ and $\Sigma_c^m$ are non-compact. This is due to collisions of the satellite with the earth, respectively the moon. However, two-body collisions can always be regularized. We explain here the regularization technique of Moser, in which one interchanges the roles of position and momentum. For any base point $b \in \mathbb{R}^2$ We consider the symplectomorphism
	$$\sigma_b \colon T^* \mathbb{R}^2 \to T^* \mathbb{R}^2, \quad
	(q,p) \mapsto (p,b-q).$$
	We think of $\mathbb{R}^2$ as a chart of the two-dimensional sphere $S^2$ using stereographic projection at the north pole. This gives rise to an embedding
	$$\iota \colon \mathbb{R}^2 \to S^2=\mathbb{R}^2 \cup \{\infty\}$$
	which induces a symplectic embedding
	$$\iota_* \colon T^* \mathbb{R}^2 \to T^* S^2.$$
	We consider the embeddings
	$$\iota_* \sigma_e (\Sigma^e_c) \subset T^*S^2, \quad
	\iota_* \sigma_m (\Sigma^m_c) \subset T^*S^2.$$
	Since the original momentum $p$ was always finite, both of the embeddings avoid the fiber over the north pole $T^*_\infty S^2$ which corresponds to infinite momentum. To regularize collisions we consider the closures
	$$\overline{\Sigma}_c^e:=\overline{\iota_* \sigma_e (\Sigma^e_c)} \subset T^*S^2, \quad
	\overline{\Sigma}_c^m:=\overline{\iota_* \sigma_m (\Sigma^m_c)} \subset T^*S^2.$$
	\begin{thm}[\cite{AFKP12}]
		The regularized energy hypersurfaces $\overline{\Sigma}_c^e$ and $\overline{\Sigma}_c^m$ are fiberwise starshaped.	
	\end{thm}
	This theorem implies that the regularized hypersurfaces are diffeomorphic to the unit tangent bundle of the two-dimensional sphere which itself is diffeomorphic to three-dimensional real projective space $\mathbb{R}P^3$. Moreover the Liouville one-form $\lambda \in T^* S^2$ restricted to
	$\overline{\Sigma}_c^e$, respectively $\overline{\Sigma}_c^m$ gives a contact form on these spaces. The embedding into $T^* S^2$ gives a filling so that the underlying contact structures are actually tight. Since there is just one tight contact structure on $\mathbb{R}P^3$ we can refer to this contact structure as the tight contact structure on $\mathbb{R}P^3$. 
	\\ \\
	The fibers over the north pole $\mathcal{C}_c^e:=\overline{\Sigma}_c^e \cap T_\infty S^2$
	and $\mathcal{C}_c^m:=\overline{\Sigma}_c^m \cap T_\infty S^2$ correspond to collisions with the earth, respectively the moon. It holds that
	$$\overline{\Sigma}^e_c=\iota_* \sigma_e(\Sigma^e_c) \cup \mathcal{C}_c^e, \quad
	\overline{\Sigma}^m_c=\iota_* \sigma_m(\Sigma^m_c) \cup \mathcal{C}_c^m.$$
	Note that the two collision circles $\mathcal{C}_c^e$ respectively
	$\mathcal{C}_c^m$ are Legendrians. Geometrically the regularization corresponds to compactify the energy hypersurface by adding a Legendrian circle.
	\\ \\ 
	The other fibers $T^*_p S^2$ for
	$p \in \mathbb{R}^2 \subset S^2$ have from the point of physics a bit strange interpretation. They correspond to some fixed momentum $p$ but arbitrary position $q$. 
	\\ \\
	What is however, remarkable is that the antisymplectic involution
	$\rho$ extends to the regularizations. Consider the meridian
	$\{0\}\times \mathbb{R} \cup \{\infty\} \subset \mathbb{R}^2 \cup \{\infty\}=S^2$
	Let $\bar{\rho} \colon T^*S^2 \to T^* S^2$ be the antisymplectic involution which is obtained as the composition of the symplectic involution induced by reflection at this meridian and the antisymplectic involution mapping each cotangent vector to its inverse. Then the regularizations interchange the antisymplectic involution $\rho$ with the antisymplectic involution $\bar{\rho}$. Note that the extended involution leaves the collision Legendrians invariant.   
	
	\section{Contact structures for Hamiltonian manifolds}\label{sec:contact_str_ham}
	
	Before discussing contact structures for the restricted three-body problem above the first critical value, we first discuss abstractly contact structures for Hamiltonian manifolds and their behaviour under antisymplectic involutions. 
	Assume that $\Sigma$ is a closed, connected, orientable odd-dimensional manifold of dimension $2n-1$. A \emph{Hamiltonian structure} on $\Sigma$ is a two-form
	$\omega \in \Omega^2(\Sigma)$ satisfying the following two conditions
	\begin{description}
		\item[(i)] $\omega$ is closed, i.e. $d\omega=0$,
		\item[(ii)] The kernel $\mathrm{ker}(\omega)$ is a one-dimensional distribution on $T\Sigma$.
	\end{description}
	For any $x \in \Sigma$ the two-form $\omega_x$ induces a symplectic form on the vector space $T_x \Sigma/\mathrm{ker}(\omega)$. In particular, the quotient bundle $T\Sigma/\ker(\omega)$ is canonically oriented by the Hamiltonian structure. Since $\Sigma$ is orientable the line bundle $\ker(\omega)$ is orientable as well and a choice of orientation on $\Sigma$ induces an orientation on $\ker(\omega)$. We refer to the tuple $(\Sigma,\omega)$ as a 
	\emph{Hamiltonian manifold}. 
	\\ \\
	A \emph{contact form} on a Hamiltonian manifold $(\Sigma,\omega)$ is a one-form $\lambda \in \Omega^1(\Sigma)$ satisfying the following conditions
	$$d\lambda=\omega, \quad \lambda \wedge \omega^{n-1} \neq 0.$$
	These two conditions imply
	$$\lambda \wedge (d\lambda)^{n-1} \neq 0$$
	so that a contact form on the Hamiltonian manifold $(\Sigma,\omega)$ is in particular a contact form on the manifold $\Sigma$. However, on Hamiltonian manifolds a contact form additionally has to satisfy some compatibility with the Hamiltonian structure $\omega$. 
	\\ \\
	A \emph{real Hamiltonian manifold} is a triple $(\Sigma,\omega,\rho)$ where
	$(\Sigma,\omega)$ is a Hamiltonian manifold and $\rho \colon \Sigma \to \Sigma$ is a smooth involution satisfying $\rho^*\omega=-\omega$ and which reverses the orientation on the line bundle $\ker(\omega)$. Suppose now that 
	$(\Sigma,\omega,\rho)$ is a real Hamiltonian manifold and $\lambda \in \Omega^1(\Sigma)$ is a contact form for $(\Sigma,\omega)$. We define
	$$\lambda^\rho:=\frac{1}{2}(\lambda-\rho^*\lambda) \in \Omega^1(\Sigma).$$
	Note that since $\rho$ is an involution $\lambda^\rho$ is antisymmetric with respect to $\rho$, i.e.
	$$\rho^*\lambda^\rho=-\lambda^\rho.$$
	We have the following lemma, which tells us that if a real Hamiltonian manifold
	$(\Sigma,\omega,\rho)$ admits a contact form it admits as well an antisymmetric contact form defining the same contact structure on $\Sigma$.
	\begin{lemma}\label{real}
		The antisymmetric one-form $\lambda^\rho$ is as well a contact form on
		$(\Sigma,\omega)$ and the contact structures $(\Sigma,\ker \lambda)$ and
		$(\Sigma,\ker \lambda^\rho)$ agree. 
	\end{lemma}
	\textbf{Proof: } We abbreviate
	$$\mathrm{vol}:=\lambda\wedge \omega^{n-1}=\lambda \wedge (d\lambda)^{n-1} 
	\in \Omega^{2n-1}(\Sigma)$$
	the volume form on $\Sigma$ induced from $\lambda$. Since $\rho$ is orientation reversing on the line bundle $\ker(\omega)$ and satisfies 
	$\rho^*\omega=-\omega$ we see that $\rho$ is orientation preserving 
	iff $n$ is even, which we write as
	$$(-1)^n \rho^* \mathrm{vol}>0.$$
	We have
	$$-\rho^* \lambda \wedge \omega^{n-1}=(-1)^n \rho^* \lambda \wedge
	(\rho^* \omega)^{n-1}=(-1)^n \rho^*(\lambda \wedge \omega^{n-1})=
	(-1)^n \rho^* \mathrm{vol}>0.$$
	For $t \in [0,1]$ we abbreviate
	$$\lambda_t:=(1-t)\lambda-t\rho^*\lambda.$$
	The above computation shows that
	$$\lambda_t \wedge \omega^{n-1}=(1-t) \mathrm{vol}+t(-1)^n \rho^*\mathrm{vol}>0.$$
	Moreover, we have
	$$d\lambda_t=(1-t)d\lambda-t \rho^*d\lambda=(1-t)\omega-t\rho^*\omega=
	(1-t)\omega+t\omega=\omega$$
	so that for every $t \in [0,1]$ is a contact form for $(\Sigma,\omega)$. 
	In particular,
	$$\lambda^\rho=\lambda_{1/2}$$
	is a contact form on $(\Sigma,\omega)$ and since the family $\lambda_t$ gives rise to a smooth homotopy between $\lambda$ and $\lambda^\rho$ it follows from Gray stability that the contact structures $(\Sigma,\ker\lambda)$ and
	$(\Sigma,\ker \lambda^\rho)$ agree. This finishes the proof of the lemma. \hfill $\square$
	\\ \\
	We further note the following. If $L \subset \Sigma$ is a Legendrian submanifold for $\lambda$ which is invariant under $\rho$, then $L$ is as well Legendrian for $\lambda^\rho$.
	
	\subsection{Contact structures above the first critical value}\label{sec:above}
	
	For energies $c \in (U(\ell_1),U(\ell_2))$ between the first and second critical value, the Hill's region has just two connected components
	$$\mathfrak{K}_c=\mathfrak{K}_c^{e,m} \cup \mathfrak{K}_c^u,$$
	where the first one $\mathfrak{K}_c^{e,m}$ is bounded and contains both the earth and the moon in its closure, while the second one is unbounded. Again, we are interested in the bounded component and abbreviate for the connected component of the energy hypersurface lying above the bounded component 
	$$\Sigma_c^{e,m}:=\{(q,p) \in \Sigma_c: q \in \mathfrak{K}_c^{e,m}\}.$$
	The Lagrange point $\ell_1$ is a saddle point of the effective potential and therefore its lift $\Pi^{-1} \ell_1$ a critical point of Morse index one of the Hamiltonian $H$. Therefore, crossing the first critical value means topologically taking the connected sum of the components of the energy hypersurface around the earth and the one around the moon. If one regularizes collisions with the earth and the moon, one obtains a closed energy hypersurface $\overline{\Sigma}^{e,m}_c$ which is diffeomorphic to the connected sum of two three-dimensional real projective spaces
	$$\overline{\Sigma}^{e,m}_c=\mathbb{R}P^3 \# \mathbb{R}P^3.$$
	\begin{thm}[\cite{AFKP12}]
		There exists $\epsilon \in (0,U(\ell_2)-U(\ell_1)]$ such that for every 
		$c \in (U(\ell_1), U(\ell_1)+\epsilon)$ as a Hamiltonian manifold 
		$\overline{\Sigma}^{e,m}_c$ admits a contact form $\lambda$, such that
		the contact structure $(\overline{\Sigma}^{e,m}_c, \ker \lambda)$ equals contact connected sum of two tight $\mathbb{R}P^3$.	
	\end{thm}
	The contact structure is obtained by careful interpolation in the neck region around the Lagrange point $\ell_1$ between the contact structures around the earth and the moon. In particular, the collision circles $\mathcal{C}_c^e$ and
	$\mathcal{C}_c^m$ are still Legendrian submanifolds for $\lambda$. The antisymplectic involution $\rho$ of the restricted three-body problem restricts and extends to a real structure on $\overline{\Sigma}^{e,m}_c$. By Lemma~\ref{real} we can assume without loss of generality that maybe after averaging, the contact form is antiinvariant under the involution. Since the Legendrians $\mathcal{C}_c^e$ and $\mathcal{C}_c^m$ are invariant under the involution, they stay Legendrians after averaging the contact form.  
	
	{\footnotesize
		\bibliography{citations}
		\bibliographystyle{alpha}
	}
\Address

\end{document}